\documentclass{article}
\usepackage[utf8]{inputenc}

\begin{document}

\title{Maximum Relative Divergence Principle for Grading Functions on Power Sets}

\author{Alexander Dukhovny \\
\small{Department of Mathematics, San Francisco State University} \\
\small{San Francisco, CA 94132, USA} \\
\small{dukhovny@fsu.edu} \\ }

\date{\today}          

\maketitle

\begin{abstract}

The concept of Relative Divergence of one Grading Function from another on a given set is extended here from totally ordered chains to subset inclusion-ordered power set of a finite event space. In particular, using general (non-additive) measures as grading functions, Shannon's Entropy concept is extended to such sets. Based on that, Maximum Relative Divergence Principle is introduced as a tool for determining the "most reasonable" grading function and used in applications where that function is supposed to be "element-additive" or "cardinality-dependent" under application-specific linear constraints.

\end{abstract}

\section{Introduction}
\label{sec:Intro}

In numerous fields of mathematics the Insufficient Reason Principle has been stated as the Maximum Entropy Principle (MEP). According to MEP, the "most reasonable" way (which is to say, using the fewest extra assumptions) to determine missing pieces of the needed probability distribution is to maximize, under some application-specific constraints, Shannon Entropy functional (see, e.g., \cite{Shannon}). That approach has proved effective in so many cases that references are just too many to quote.\par

There have been numerous generalizations of the original Shannon Entropy formula. The list includes relative entropy, Kullback-Leibler divergence, partition entropy,  Kolmogorov-Sinai entropy, topological entropy, entropy of general non-probabilistic measures (capacities) and great many others  (see, e.g., references [5-11] and a recent review in \cite{Review}.

In our preceding article \cite{DArXiv} Shannon Entropy was generalized to the concept of Relative Divergence (RD) of one Grading Function (GF) from another on a totally ordered set (chain), reducing to Shannon Entropy in a special case. (The term Relative Divergence was chosen in keeping with Kullback-Leibler Divergence - see \cite{K-L} - also known as Relative Divergence of probability measures.)

Here we begin the process of generalization of the Relative Divergence concept and the associated Maximum Relative Divergence Principle (MRDP) to partially ordered sets. Doing that, we show that 

1. both of those concepts reduce to Shannon Entropy and Maximum Entropy Principle when dealing with a probability theory problem;

2. conclusions made by using MRDP in new problems agree with "common sense" ones - where such are available;

3. MRDP can be effectively used in new applications.

We start the process by working with a power set $W = 2^X$ of the event space $X$ ordered by subset inclusion. In doing that, we will be following the way the concept of General Entropy of General (nonadditive) Measures (see \cite{DFuzzy}) was developed for power sets. 

Using a normalized General Measure $\mu (w)$ as a subset $w$ Grading Function $F(w$, it follows that its relative divergence from the "subset cardinality" grading function $|w|$ reduces to the minimum of Shannon Entropy values of all probability distributions said to be "subordinate" to that general measure (as befits a generalized concept).

We then apply Maximum Relative Divergence Principle (MRDP) to some problems outside of the probability theory - the natural domain of MEP - and show that it does lead to "most reasonable, common sense" results. In particular, MRDP will be applied to some problems arising in Operations Research, introducing a new tool for their analysis and providing some new results.

The general setup of Relative Divergence in  \cite{DArXiv} begins as follows: let $W$ be a totally ordered set and use $\prec$ to denote the ordering relation of its elements. A real-valued function $F : W \to R$ is said to be a Grading Function on $W$ if 

$w \prec v \iff F(w) < F(v)$ for all $w, v \in W $.

(In particular, when $W$ is countable the image $i(w)$ of an element $w$ under an order-preserving isomorphism $i: W \longrightarrow Z$ is said to be a "natural" GF on $W$.)

As such, the inverse function $w = F^{-1}(u)$  is defined  for all $u \in im(F) \subset R$.

When grading functions $F(w)$ and $G(w)$ are defined on $W$, relative divergence RD of $F$ from $G$ over $W$ is defined by one of the two formulas below depending on the nature of $W$ and $im(F)$.

When $W = \{ \ldots , w_{-1}, w_0 , w_1 , \ldots \}$ is countable the RD of $F$ from $G$ on $W$ is defined as

\begin{equation} \label{discrete RD}
\mathcal{D}(F \Vert G) \vert_W = \ \sum_{k= - \infty}^\infty\
\ln  \left( \frac{\Delta_k G}{\Delta_k F }\right) \Delta_k F , 
\end{equation}

\noindent where  

$\Delta_k F = F(w_k) - F(w_{k-1}), \quad \Delta_k G= G(w_k) - G(w_{k-1}), \quad  k = \ldots,-1,0,1,\ldots$,
\noindent assuming absolute convergence of the series.

When $im(F)$ is an interval, relative divergence of $F$ from $G$ on $W$ is defined as

\begin{equation}\label{continuous RD}
\mathcal{D}(F \Vert G) \vert_W =
\int_{im(F)} \ln \left( \frac{d}{du}(G(F^{-1}(u))) \right)du,
\end{equation}

\noindent assuming absolute convergence of the integral.

It follows trivially from (\ref{discrete RD}) and (\ref{continuous RD}) that $\mathcal{D}(F \Vert G)$ does not change if either (or both) of the grading functions is shifted by a constant. As such, when $W$ is well-ordered (that is, it has a minimal element $m_W$) then (if needed) with no generality loss it can be assumed that $F(m_W) = G(m_W) = 0$.

The following special case of (\ref{discrete RD}) highlights the connection between Relative Divergence and Shannon Entropy:

When $W = \{w_0, w_1, \ldots, w_n \}$ and $I(w_i) = i, \quad i=1, 2, \ldots, N$ (that is, $I$ is the ordinal function on $W$), and 
$F(w_0)=0$, \quad $F(w_n)=1$,  then

\begin{equation} \label{ordinal RD}
\mathcal{D}(F \Vert I) \vert_W = - \sum_{k=0}^n \ln  (\Delta_k F ) \Delta_k F , 
\end{equation}

\noindent and is called Shannon Entropy $\mathcal{H}(F)$ of $F$ on $W$.

Taking the limit as $n \longrightarrow \infty$, when $F(w_n) \longrightarrow 1$, formula (\ref{ordinal RD}) extends to the case of a well-ordered countable $W$ (assuming that the limit exists).

When $F$ is a bounded grading function on a countable $W$, we denote

$m_F = \inf{im(F)}, \quad M_F = \sup{im(F)}, \quad M_F - m_F = \Delta _W F$.

\noindent and introduce the "normalized" grading function $\hat{F}$ by 

$\hat{F} (w)= \frac{F(w)-m_F}{\Delta_W F}$, 

\noindent which can be also expressed as

$F(w) = (\Delta_W F) \hat{F}(w) + m_F$ \quad or  \quad
$\Delta_k F = (\Delta_W F )\Delta_k\hat{F}$.

As such, for the special case where $G = I$, it follows directly from (\ref{ordinal RD}) that

\begin{equation} \label{ordinal entropy}
\mathcal{D}(F \Vert I) \vert_W = (\Delta_W F) \mathcal{H}(\hat{F}) + (\Delta_W F) \ln{(\Delta_W F)}.
\end{equation}

Also, equation (\ref{discrete RD}) yields the following formula that connects the RD of two grading functions to the RD of their normalized versions:

\begin{equation} \label{normalized RD}
\mathcal{D}(F \Vert G) \vert_W = 
-(\Delta_W F) \ln  \left( \frac{\Delta_W F}{\Delta_W G }\right)   +  
(\Delta_W F) \mathcal{D}(\hat{F} \Vert \hat{G})
\end{equation}

Extending formula (\ref{discrete RD}) to partially ordered sets (posets), our approach here is to look at posets as unions of (possibly overlapping) maximal linear orders (chains). We concentrate on the cases where the poset $W$ is countable, well-ordered and connected. Namely, we do it here for power sets of subsets ("events") of an "event" space of elements ("outcomes") ordered by subset inclusion.

In Section 2, when a general (non-additive) measure is defined on $W$, it is treated as a subset-grading function. Based on that, we extend Shannon Entropy to such cases and develop Maximum Relative Divergence Principle (MRDP) as a tool to make "most reasonable" decisions on constructing such measures.

In section \ref{MRDP for C-D case} MRDP is applied to cases where the "most  reasonable" subset-grading function is supposed to have values dependent only on the subset's cardinality with some of those values pre-specified by the nature of the application.

In section \ref{MRDP for E-A under cons}, we specify MRDP to applications where the "most reasonable" grading function must be "element-additive" and should satisfy some linear constraints on its increments.

In particular, in that section, in the context of "resource distribution" applications, we use the the results of the section for cases where the "most reasonable" (element-additive) grading function must have pre-specified values ("resource quotas" on subsets partitioning the event space. Furthermore, in section \ref{quotas and costs}, we explore the case where that grading function must also satisfy the fixed "subset cost" constraints.

\section{Relative Divergence on Power Sets}\label{RD on Psets}

 The main object of interest in this section is the concept of Relative Divergence of grading functions on the power set $W = 2^X$, where  $X = \{x_i, \quad i = 1, 2, \ldots, |X|\}$ is the event space. As usual, $W$ is ordered by subset inclusion, so it possesses the least element $w_0 = \emptyset$ and each element of $W$ belongs to at least one maximal chain ($MC$) in $W$ (which must start at $w_0$). 
 
 Our approach here is to treat $W$ as the union of all its maximal chains:

$W=\cup \{\forall MC \subset W \} $.

Accordingly, when grading functions $F$ and $G$ are defined on $W$, we propose  to construct 
$\mathcal{D}(F \Vert G) ) \vert _W$, 
their relative divergence on the entire $W$ as an aggregate of relative divergences of those grading functions 
$\mathcal{D}(F \Vert G) \vert_{MC}$ 
along all maximal chains $MC$ in $W$.

When $W$ is a general poset, the structure of the order  relation on $W$ and interdependence of the grading functions may complicate aggregation or even make it infeasible.

Here, where $W$ is a power set, its special properties suggest a meaningful way for such aggregation: all maximal chains are specified by permutations of $X$, they start at the same  point $w_0 = \emptyset$,  have the same number of elements $n+1$ and end in the same element $w_n = X$, so the grade spreads of the GFs along each maximal chain are all equal:

$\Delta_{MC} F = \Delta_W F, \quad   \forall{MC} \subset W$.

The aggregation method chosen here follows the one proposed in \cite{DFuzzy} to facilitate Maximum Entropy Principle for a general measure $\mu$ on $W$: its Shannon entropy $\mathcal{H} (\mu) \vert_W$ should be taken as the minimum over all maximal chains $MC$ in $W$ of Shannon Entropy values of probabilistic measures $\mu_{MC}^s$ (said to be "subordinate to $\mu$ on $MC$"). Those measures are completely and uniquely determined by having their values equal to the values of $\mu$ on $MC$. 

Now, formula (\ref{ordinal entropy}) relates the relative divergence $F \Vert I \vert_W$ to Shannon entropy of the general measure $\mu$ generated on $W$ by the normalized  $\hat{F}$. In line with the method of \cite{DFuzzy}, in this paper we therefore propose the following definition:

\begin{equation} \label{RD on W}
\mathcal{D}(F\Vert G)\vert_W = 
\min_{MC\subset W} \mathcal{D}(F\Vert G)\vert_{MC}
\end{equation}

To facilitate that definition, we represent each maximal chain $MC$ in $W$ as a sequence $\{w_i, \quad i=0, 1, \ldots, n \}$ determined by a permutation $ \{ x_{k(i)} , \quad i = 1, \ldots n \} $ of the elements of the event space $X$ as follows: 

$w_0 = \emptyset, \quad w_i = w_{i-1} \cup x_{k(i)},  \quad i = 1, \ldots n.$

Accordingly, we define 

$f (x_{k(i)}) = F(w_i) - F(w_{i-1}), \quad i = 1, \ldots, n,$ 

\noindent the "increment" function of $F$ along the chain $MC$. As a grading function, $F$ is set-monotonic, so its increment function assumes only positive values. $MC$ being a maximal chain, by its definition, $f$ is defined at all elements of $X$. 

Now we can define another grading function on $W$:

$F^s _{MC} (w) = \sum_{ x \in w} f (x), \quad \forall w \in W.$

That relation makes $F^s _{MC}$ an "element-additive" grading function. Its values are uniquely determined by the values of $F$ on the chain $MC$. Moreover,

$F^s _{MC} (w) = F(w), \quad \forall w \in MC$, 

\noindent so $F^s _{MC}$ is said to be "subordinate to $F$ on $MC$ ".

When $W$ is a power set of $X$, a grading function $F$ is said to be cardinality-dependent if 
$F(w) = |w|, \quad \forall w \in W)$.

In particular, the cardinality of a subset $N(w) =|w|, \quad \forall w\in W$. is itself a "natural" grading function on $W$.

Describing an MC in $W$ as a sequence $\{ w_i, \quad i=0, 1, \ldots, n \}$ where $N(w_i) = i$, it follows from (\ref{ordinal entropy}) and (\ref{RD on W}) that

\begin{equation} \label{RD and H}
\mathcal{D}(F \Vert N) \vert_W = \mathcal{H}(\hat{F}) \vert_W +
\Delta_W F \ln{\Delta_W F}.
\end{equation},

\noindent where
$\mathcal{H}(\hat{F}) \vert_W = 
\min_{MC\subset W} \mathcal{H}(\hat{F})\vert_{MC}$

Formula (7) generalizes formula (4) for a linearly ordered $W$ to the case where $W$ is a power set of a finite $X$. Similar to \cite{DFuzzy}, we define a grading function $F$ to be "equilateral" when 
$\mathcal{D}(F\Vert N)\vert_{MC}$ is the same for all maximal chains $MC \in W$.

It follows directly from the definitions of element-additive and cardinality-dependent grading functions that they are equilateral. As well, any linear transformation of any equilateral grading function $F$ (say, the corresponding normalized grading function $\hat{F}$) is also equilateral.

As such, it opens a way to generalize Maximum Entropy Principle (MEP) for probability distributions to Maximum Relative Divergence Principle (MRDP) for grading functions:

MRDP: An "admissible" (satisfying the constraints of the problem) grading function $F$ on a power set $W$ is said to be "the most reasonable" (within the fixed grading interval) if it gives maximum to 
$\mathcal{D}(F \Vert N) \vert_ W $.

Now, using MRDP involves analysis of all 
$\mathcal{D}(F \Vert N) \vert_ MC, \forall{MC \in W} $.

In general, using MRDP (same as MEP), leads to a nonlinear problem on the constraint-imposed domain. The added complexity of that problem for MRDP is that it involves considering all maximal chains in the event space. 

In this paper we explore applications where the emerging grading functions are equilateral, so only one maximal chain needs to be analyzed. That can often be done using techniques similar to the ones generally used to treat problems related to Shannon Entropy theory (see, e.g., \cite{ElementsIT}). In particular, of special significance is the following fundamental result (see, e.g.,  \cite{Jaynes}).

Lemma 1. For a probability distribution 

$\{p_i\},0\leq p_i\leq 1, \quad \sum_{i=1}^n p_i = 1, \quad i = 1, \ldots, n$,

\noindent the maximum value of Shannon entropy of that distribution

$\mathcal{H}\quad  = -\sum_{i=1}^n p_i\ln{p_i} \quad = \ln{n}$

\noindent is attained when $p_i = \frac{1}{n}, \quad i = 1, \ldots, n$.

\section{MRDP for cardinality - dependent subset-grading function}
\label{MRDP for C-D case}

Using MRDP to find the "most reasonable" grading function $F(w)$ on $W = 2^X$ leads to a direct formula where

$F(w) = F(|w|, \quad \forall w \in W, \quad F(0) = 0, \quad F(|W|) = M$.

A typical application may arise in the context of the queuing theory with group service. Say, testing a sample group $w$ of subjects from a population of size $N$ for the presence of an "attribute-positive" subject, when the "cost" $F(w)$ of the group test  depends only on the size of the group. 

In particular, let us apply MRDP in the case where some values of $F(|w|)$ are pre-specified, that is,

$F(n_k)=M_k, \quad k=1, \ldots, K$, 

\noindent where we define 

$n_0 = 0, \quad M_0 = m,\quad n_K = n = |X|, \quad M_n = M$.

Denote index intervals 
$I_k = ( n_{k-1}, n_k ], \quad k = 1, \ldots, K$.

\noindent and introduce quantities

$q_{i,k} = \frac{F(i) - F(i-1)}{ M_k - M_{k-1}}, \quad i \in I_k$.

For a cardinality-dependent $F(w)$ the $\mathcal{D}(F \Vert N) \vert_{MC}$ is the same for each maximal chain $MC$ (since successive subsets in any $MC$ differ by just one element), so the MRDP problem presents as follows: 

Find the $q_{i,k} \geq{0}, \quad i \in I_k, \quad k = 1, \ldots, K$

\noindent to maximize $ \quad -\sum_{k=1}^K \sum_{i \in I_k} q_{i,k} \ln (q_{i,k} ) $

\noindent subject to $\quad \quad \sum_{i \in I_k} q_{i,k} = 1, \quad k = 1, \ldots,K$.

The additive form of the maximized expression leads to $K$ independent maximization problems for each $k = 1, \ldots, K$. Using Lemma 1, the (unique) solution of each one of them is

$q_{i,k} = \frac{1}{n_k - n_{k-1}}, i \in I_k$.

Respectively, the overall solution of the MRDP problem presents as a piece-wise linear function as follows:

\begin{equation} \label{C-D solution)}
F(w)=a_k+b_k |w|, \quad |w| \in I_k, \quad k=1, \ldots, K, \quad \forall w \in W
\end{equation}

\noindent where
$\quad b_k= \frac{ M_k - M_{k-1} }{n_k - n_{k-1}}, \quad a_k = M_k - b_k n_{k-1}.$ 

In particular, when only the total testing cost $M$ of the entire population is specified,  $K=1, m = 0$, and formula (\ref{C-D solution)}) reduces to

$F(w) = M \frac{|w|}{n},\quad \forall w \in W$.

\section{MRDP for element-additive grading functions under linear constraints on the increments}
\label{MRDP for E-A under cons}

Here we consider the case where the subset-grading function must be element-additive:

$F(w) = \sum_{x_j \in w} f(j)$, 
where $f(j) \geq 0, \quad F(w_0) = 0, \quad F(X) = M$

\noindent and, in addition to 

$f(1) + \ldots + f(n) = M, \quad n = |X|$

\noindent another set of $K$ linear constraints must be satisfied:

$\sum_{j=1}^n a_{k,j} f(j) = M_k, \quad k = 1, \ldots, K$.

As such, the MRDP problem for this case look as follows:

Find the values of $f_i, \quad i = 1, \ldots, n=|X|$, to maximize 

\begin{equation}\label{DFN+}
\mathcal{D}(F \Vert N) \vert_ W =-\sum_{k=1}^K \sum_{i \in I_k} f_i \ln (f_i) 
\end{equation}

subject to 

\begin{equation}\label{constraints}
A \vec{f} = \vec{M},
\end{equation}

\noindent where $A$ is the matrix of constraint coefficients whose $i$-th column is

$\vec{A}_i =  [a_{1,i}, \ldots, a_{k,i}]^T, \quad i = 1, \ldots, n$,
\quad and

$\vec{f} = [f(1), \ldots, f(n)]^T, \quad \vec{M} = [M_1, \ldots, M_K]^T $

Using Lagrange multipliers method to find a solution at an interior point of the constraints-imposed domain and combining the multipliers for each constraint into a column-vector 
$\vec{\lambda} = [\lambda_1, \ldots, \lambda_K]$, it follows that

$-\ln (f_i) - 1 - \vec{\lambda} \vec{A}_i = 0, \quad i = 1, \ldots, n$,

\noindent so that

\begin{equation}\label{f-formula}
\vec{f} = e^{-1}[e^{- \vec{\lambda} \vec{A}_i}, i = 1, \ldots, n]^T,
\end{equation}

\noindent where $\vec{\lambda}$ is to be found from the system of equations that arises from (\ref{constraints}):

\begin{equation}\label{lambdas}
A[e^{- \vec{\lambda} \vec{A}_i}, \quad i = 1, \ldots, n]^T = e\vec{M}
\end{equation}

The way to find a solution of the system for Lagrange multipliers which would, when used in (\ref{f-formula}), yield $\vec{f}$ as an interior point of the constraints-imposed domain depends on the structure of the constraints. That, in turn, comes from the nature of the application. Here we present several applications where the constraints allow for analytic solutions.

Consider, for example, the problem where the "most reasonable" distribution of some "resource" between the "users" in $X$ must be found under "pre-fixed" group quotas". Denoting $f(x_j)$ the resource allocation to user $x_j \in X, \quad j = 1, \ldots, n$, for a user subset $w \in W$ its subset allocation, if defined as a sum of the individual $w$-members' allocations,

$F(w) = \sum_{x_j \in w} f(x_j)$,

\noindent which makes it an element-additive grading function on $W$.

We further assume that the "users" population $X$ is partitioned into disjoint subsets $X = X_1 \cup, \ldots,  \cup X_k, \quad k = 1, \ldots, K$, whose subset resource quotas ${M_1, \ldots, M_K}$ are fixed.

Under those assumptions, denoting

$p_j = \frac{f(x_j)}{M_k}, \quad \forall x_j \in X_k, \quad k=1, \ldots K$,

\noindent the MRDP problem of the case  presents as follows: maximize

\begin{equation} \label{quotas}
\mathcal{D}(F \Vert N) \vert_ W =
\sum_{k=1}^K (M_k\ln{M_k - \sum_{x_j \in X_k} p_j\ln{p_j} )}.
\end{equation},

subject to

\begin{equation} \label{quo-cons}
\sum_{x_j \in X_k} p_j = 1, \quad k = 1, \ldots, K.
\end{equation},

As such, due to the additive form of the RHS of (\ref{quotas}), the maximization problem of $\mathcal{D}(F \Vert N) \vert_ W$ breaks up into to $K$ separate independent maximization problems for each $k = 1, \ldots, K$. Using Lemma 1, the overall solution presents as the "partition-wise" uniform distribution:

\begin{equation}\label {quo-sol}
f(x_j)=\frac{M_k}{|X_k|}, \quad \forall j: x_j \in X_k, \quad k = 1, \ldots, K.
\end{equation}

(Also, in the framework of the "refinement" procedure, when an element-additive grading function $F$ is defined on $2^X, X =(X_1, \ldots, X_K)$, and $F(X_k) = M_k, \quad k = 1, \ldots, K$, and each of the elements of $X$ is becomes a collection of its own "smaller" "sub-elements", the new power set $W$ emerges. Looking for the "most reasonable" way to expand $F$ on $W$ as an element-additive grading function preserving the "old" values, the same MRDP problem arises.)

In particular, when only the overall population resource allocation $M$ is specified, in agreement with the famous probability theory result when $M=1$ (see, e.g., \cite{Jaynes} ),

$f(x_j) = \frac{M}{|X|}, \quad j = 1, \ldots, |X|$,

\section{MRDP for resource distribution with quotas and "costs"}
\label{quotas and costs}

Using MRDP to decide on the "most natural" grading function may require  including additional constraints on its values as needed to reflect the nature of the application. Here we consider the case with linear constraints which can be interpreted as "cost" requirements.

In the context of the model in section \ref{MRDP for E-A under cons}, consider the case where the values of the grading function for a set of disjoint subsets partitioning the "users" population $X$ into disjoint subsets $X = X_1 \cup, \ldots,  \cup X_k, \quad k = 1, \ldots, K$. This time, not only subset resource quotas ${F(X1) = M_1, \ldots, F(X_K) = M_K}$ but also the total "costs" $Q_k$ of those quotas  modeled as sums of individual element $x_i$ costs proportional to the allocated element resource with the rate $r_i$. As such, the MDRP problem for the case look as follows:

Find the values of $f_i, \quad i = 1, \ldots, |X|$ to maximize 

\begin{equation}\label{MRDP for quotas and costs}
\mathcal{D}(F \Vert N) \vert_ W =-\sum_{k=1}^K \sum_{i \in I_k} f_i \ln{f_i} 
\end{equation}

\noindent subject to 

\begin{equation}\label{totals}
\sum_{x_i \in X_k} f_i = M_k, \quad k = 1, \ldots, K.
\end{equation}

\noindent and

\begin{equation}\label{costs}
\sum_{x_i \in X_k} r_i f_i = Q_k, \quad k = 1, \ldots, K.
\end{equation}

The additive form of the RHS of equation (\ref{MRDP for quotas and costs}) and non-overlapping structure of the constraints for each subset $X_k$ break down the entire problem into separate and independent maximization problems (\ref{MRDP for quotas and costs}), (\ref{totals}), (\ref{costs}) for each partition subset 
$X_k, \quad k=1, \ldots, K$. 

As such, we first approach the case where $K=1$ and drop index $k$. To find a possible interior maximum in the constraint-imposed domain, using Lagrange multipliers' method and denoting by $\lambda$ and  $\alpha$ Lagrange multipliers for equations (\ref{totals}), and (\ref{costs}), it follows that

$-\ln{f_i} - 1 - \lambda - r_i\alpha  = 0, \quad \forall{i: x_i \in X}$

\noindent which results in

\begin{equation}\label{TC1 solution}
f_i = e^{ -1 - \lambda - r_i \alpha}, \quad i = 1,\ldots, n =|X|
\end{equation}

\noindent with $\lambda, \alpha$ to be determined from the system of equations arising from equations (\ref{totals}) and (\ref{costs}) upon using formula (\ref{TC1 solution}):

$\sum_{i=1}^n e^{ -1 - \lambda - r_i \alpha} = M$,

$\sum_{i=1}^n r_i e^{ -1 - \lambda -r_i \alpha} = Q$.

From those equations, upon some trivial algebraic transformations and denoting $r_0 = \frac{Q}{M}$ - "the constraints-imposed cost rate", and $ d_i = r_i - r_0$ - deviations of individual elements' cost rates from $r_0$, we obtain an equation to find $\alpha$:

$E(\alpha) = \sum_{i=1}^n d_i e^{-d_i \alpha} = 0$.

Since $E^{'} (\alpha) = - \sum_{i=1}^n (d_i)^2 e^{-d_i \alpha} \leq 0 $
it follows that: 

1. When $d_i = 0$ (that is, $r_i = r_0$), $\quad \forall i$, equation $E(\alpha)=0$ has infinitely many solutions. However, it simply means that the "total costs" constraint is satisfied by the "total $M$" constraint, so (\ref{TC1 solution}) is given by 

$f_i = \frac{M}{n},  \quad i = 1,\ldots, n =|X|$;

2. When all $d_i$, are of the same sign (but not all zeroes), equation $E(\alpha)=0$ has no solutions. It means that the resource distribution problem itself is inconsistent: the  individual elements' resource allowances and costs cannot be reconciled with the imposed total resource and cost amounts $M$ and $Q$;

3. When $d_i, \quad i = 1,\ldots\, n$ are not of the same sign equation $E(\alpha)=0$ has a unique solution (to be determined by numerical methods).

(In particular, a nontrivial tractable result emerges when the average cost rate $\bar{r} = \frac{\sum_{i=1}^{n}r_i}{n} = r_0$, 
the constraint-imposed cost rate. Then the unique solution of $E(\alpha)=0$ is $\alpha = 0$, so once again
$f_i = \frac{M}{n},  \quad i = 1,\ldots, n$.)

\section{Conclusion}
\label{end}

In summary: the concept of Relative Divergence (RD) of grading functions introduced in \cite{DArXiv} for totally ordered sets is extended here to power sets of finite event spaces as a generalization of Shannon Entropy for probability distributions on event spaces. 

Based on that, Maximum Relative Divergence Principle (MRDP) was introduced, studied and applied here to related Operations Research problems as an extension of Maximum Entropy Principle (MEP) for Probability Theory problems. It was shown that: 

1. MRDP reduces directly to MEP when the problem in question is similar to a Probability Theory problem;

2. MRDP applied to some test problems from Operations Research led to "most reasonable", "common sense" conclusions;

3. MRDP opens another way to analyze some Operations Research applications (such as group service in Queuing Theory, resource distribution under constraints, etc.).

\end{document}